\documentclass[10pt]{article}
\usepackage{amsmath}
\usepackage{amsfonts}
\usepackage{amssymb}
\usepackage{comment}
\usepackage{graphicx,color}
\long\def\proof#1{\removelastskip\vskip\baselineskip\relax\noindent{\it
Proof\if!#1!\else\ \ignorespaces#1\fi.\ }\ignorespaces}

\newcommand{\lgs}[2]{\mbox{$\left(\frac{#1}{#2}\right)$}}

\newcommand{\ov}{\overline}
\newcommand{\Q}{{\mathbb Q}}
\newcommand{\Z}{{\mathbb Z}}

\newcommand{\G}{\Gamma}

\newcommand{\al}{\alpha}
\newcommand{\be}{\beta}

\DeclareMathOperator{\CS}{CS}

\DeclareMathOperator{\SL}{SL}

\DeclareMathOperator{\lcm}{lcm}

\newcommand{\Proof}{{\it Proof. \/}}
\newcommand{\squareforqed}{\hbox{\rlap{$\sqcap$}$\sqcup$}}
\newcommand{\qed}{\ifmmode\squareforqed\else{\unskip\nobreak\hfil
\penalty50\hskip1em\null\nobreak\hfil\squareforqed
\parfillskip=0pt\finalhyphendemerits=0\endgraf}\fi}

\newcommand{\fp}{\qed\removelastskip\vskip\baselineskip\relax}

\newtheorem{theorem}{Theorem}[section]
\newtheorem{corollary}[theorem]{Corollary}
\newtheorem{proposition}[theorem]{Proposition}
\newtheorem{lemma}[theorem]{Lemma}
\newtheorem{definition}[theorem]{Definition}
\newtheorem{conjecture}[theorem]{Conjecture}

\newcommand{\litem}{\par\noindent\dimen0=\parindent%
    \advance\dimen0 by-4pt
               \hangindent=\dimen0\ltextindent}

\newcommand{\ltextindent}[1]{\hbox to \hangindent{#1\hss}\ignorespaces}
\newcommand{\ltextjndent}[1]{\hbox to \hangindent{#1\hss}\ignorespaces\kern-1ex}

\renewcommand{\pmod}[1]{\allowbreak\ ({\rm{mod}}\,\,#1)}

\begin{document}
\pagestyle{plain}
\title{On the Argument of the Lerch, Chowla--Selberg\\
  Formula and CM Values of $\eta(\tau)$}

\author{Henri Cohen}

\maketitle

\begin{abstract}
  We give a completely explicit Lerch, Chowla--Selberg formula without
  absolute values of the Dedekind eta function, and give some applications,
  in particular to the argument of individual CM values of the eta function.
  We also give precise conjectures on these CM values and a large table
  supporting them.
\end{abstract}  
  
\section{Introduction}

\begin{definition}
  Let $D<0$ be congruent to $0$ or $1$ modulo $4$.
  \begin{enumerate}
  \item We denote by $Q(D)$ the set of quadratic numbers
    $\tau=(-b+\sqrt{D})/(2a)$ (or, equivalently, of binary quadratic forms
    $(a,b,c):=ax^2+bxy+cy^2$) such that $b^2-4ac=D$, $a>0$, and $\gcd(a,b,c)=1$
    modulo the natural action of $\SL_2(\Z)$. A representative of an element
    of $Q(D)$ will be called a \emph{CM point} of discriminant $D$.
  \item We let as usual $h(D)=|Q(D)|$, $w(D)$ the cardinality of the
    number of roots of unity in the quadratic order of discriminant $D$
    (so $w(-3)=6$, $w(-4)=4$, and $w(D)=2$ for $D<-4$), and
    $h'(D)=h(D)/(w(D)/2)$, so that $h'(D)=h(D)$ for $D<-4$.
  \item We denote by $F(D)$ the standard set of representatives of $Q(D)$,
    in other words satisfying the additional conditions $|b|\le a\le c$ and
    $b\ge0$ if either $|b|=a$ or $a=c$.
  \end{enumerate}
\end{definition}

\begin{definition} For $D$ a negative fundamental discriminant, we define the
  Chowla--Selberg gamma quotient by
  $$\CS(D)=\left(\prod_{1\le j\le|D|}\G(j/|D|)^{\lgs{D}{j}}\right)^{1/h'(D)}\;.$$
\end{definition}

Thus, $\CS(-3)=(\G(1/3)/\G(2/3))^3$, $\CS(-4)=(\G(1/4)/\G(3/4))^2$, and
otherwise the exponent is $1/h(D)$.

The Lerch, Chowla--Selberg formula, as generalized in \cite{Nak-Tag}
is as follows:

\begin{theorem}\label{thm1} Let $D<0$ be congruent to $0$ or $1$ modulo $4$,
and write $D=D_0f^2$ with $D_0$ a fundamental discriminant. Then
$$\left(\prod_{\tau\in Q(D)}\Im(\tau)|\eta(\tau)|^4\right)^{1/h(D)}=c(D_0,f)\dfrac{\CS(D_0)}{4\pi|D|^{1/2}}\;,$$
with
$$c(D_0,f)=\prod_{p\mid f}p^{e(p)},\text{\quad with\quad}e(p)=\dfrac{(1-p^{-v_p(f)})\left(1-\lgs{D_0}{p}\right)}{(1-1/p)\left(p-\lgs{D_0}{p}\right)}\;.$$
\end{theorem}

The original formula was for $D=D_0$ fundamental, hence $c(D_0,f)=1$. Note
a misprint in Section 10.5 of \cite{Coh}, where the factor $p-\lgs{D_0}{p}$
is incorrectly written as $p-\lgs{D_0}{p}/p$.

\smallskip

The first goal of this paper is to give a formula for the same product, but
without the modulus in front of $\eta(\tau)$. We will then give results and
precise conjectures on \emph{individual} values of $\eta(\tau)$.

Note first that by modularity,
the expression $\Im(\tau)|\eta(\tau)|^4$ is \emph{invariant} under the action
of $\SL_2(\Z)$, so taking the product over $Q(D)$ makes sense. On the
contrary, $\Im(\tau)\eta(\tau)^4$ is \emph{not} invariant, so we have to
specify which set of representatives we choose. Thus, we define the
function that we will study in the first part of this paper as follows:

\begin{definition}\label{def:pd}
  For $D<0$ congruent to $0$ or $1$ modulo $4$, we define
  $$P(D)=\prod_{\tau\in F(D)}\dfrac{\eta(\tau)^4}{|\eta(\tau)|^4}=\dfrac{\prod_{\tau\in F(D)}\Im(\tau)\eta(\tau)^4}{(c(D_0,f)\CS(D_0)/(4\pi|D|^{1/2}))^{h(D)}}\;,$$
  where $D=D_0f^2$ with $D_0$ fundamental and $F(D)$ is the standard
  set of representatives of $Q(D)$ as defined above.
\end{definition}

Theorem \ref{thm1} tells us that the two definitions of $P(D)$ are
equivalent, and evidently its modulus is equal to
$1$, so our goal is to understand its argument.

\section{The Main Theorem}

\begin{theorem}\label{thm2} Let $D<0$ be congruent to $0$ or $1$ modulo $4$.
  Define $U(D)=\emptyset$ if $v_2(D)=3$ or $v_2(D)=4$ or $D\equiv4\pmod{16}$,
  and otherwise $U(D)$ to be the set of positive divisors $d$ of $|D|$ such
  that $d\le\sqrt{|D|}$ and $\gcd(d,|D|/d)=2^v$ with $v=0$ if $D\equiv1\pmod4$,
  $v=1$ if $D\equiv12\pmod{16}$, and $v=2$ if $D\equiv0\pmod{32}$.
  
  Let $U^-(D)$ (resp $U^+(D)$) be the set of elements $d$ of $U(D)$ such that
  $d\le\sqrt{|D|/3}$ (resp $d>\sqrt{|D|/3}$). We have
  $$P(D)=e^{-(i\pi/6)(|U^-(D)|-3|U^+(D)|)}\prod_{d\in U^+(D)}\dfrac{d-\sqrt{D}}{d+\sqrt{D}}\;.$$
\end{theorem}

Warning: the stated congruences are for $D$ itself, \emph{not} for $|D|=-D$.

\Proof It is immediate to check that $\eta(-\ov{\tau})=\ov{\eta(\tau)}$.
If follows that if both $(-b+\sqrt{D})/(2a)$ and $(b+\sqrt{D})/(2a)$
belong to $F(D)$, the product of their eta values is a positive real number,
and $\eta(\tau)$ is also positive real if $b=0$. Since we are only concerned
with the argument of $P(D)$, it follows that the only $(a,b,c)$ which can
contribute to the argument are the boundary cases for which $a=c$ or $a=b$.

\smallskip

The case $a=c$ occurs for $D=b^2-4a^2=-(2a-b)(2a+b)$ and $0\le b\le a$.
Set $d=2a-b$ and $d'=2a+b=|D|/d$. The condition $b\ge0$ means that 
$d\le\sqrt{|D|}$. The condition $b\le a$ is easily seen to be equivalent to
$d\ge\sqrt{|D|/3}$. We have $a=(d'+d)/4$, $b=(d'-d)/2$, and the last condition
is that we need $a$ and $b$ to be integral with $\gcd(a,b,c)=\gcd(a,b)=1$,
in other words $\gcd(d'+d,2(d'-d))=4$.

If $D\equiv1\pmod4$, we have $d\mid D$ odd, and since $dd'=-D\equiv3\pmod4$,
we have $d\equiv1\pmod4$ and $d'\equiv3\pmod4$ or the reverse, so
$4\mid(d'+d)$ and $2(d'-d)\equiv4\pmod8$, hence the $2$-part of the gcd
is indeed equal to $4$, so the condition is $\gcd(d,d')=1$.

If $D\equiv0\pmod4$, the integrality of $b$ implies that $d$ and $d'$ must
both be even. Thus, if $D\equiv4\pmod{16}$, i.e, $|D|/4\equiv3\pmod4$,
we have $d/2\equiv1\pmod4$ and $d'/2\equiv3\pmod4$ or the reverse,
so $8\mid d'+d$ and $4\mid d'-d$, hence $8\mid\gcd(d'+d,2(d'-d))$, so the
gcd condition cannot be satisfied, in other words $U(D)=\emptyset$.

If $D\equiv12\pmod{16}$, we now have $d/2\equiv d'/2\equiv\pm1\pmod4$,
so $d'+d\equiv\pm4\pmod8$ and $2(d'-d)\equiv0\pmod{16}$, so the $2$-part
of the gcd is equal to $4$, so the condition is $\gcd(d/2,d'/2)=1$, i.e.,
$\gcd(d,d')=2$.

If $v_2(D)=v\ge3$, write $d=2^wm$ with $m$ odd, so $d'=2^{v-w}m'$ with
$m'$ odd. Exchanging $d$ and $d'$ if necessary we may assume $w\le v/2$.
If $w<v/2$, we have $v_2(d'+d)=w$ and $v_2(2(d'-d))=w+1$, so the gcd condition
implies that $w=2$, so $v\ge5$. If $w=v/2$, we have $v_2(d'+d)\ge w+1$ and
$v_2(2(d'-d))\ge w+2$, so the gcd condition implies $w=1$ hence $v=2$,
which is excluded. Summarizing, this shows that if $v=3$ or $v=4$ the
gcd condition cannot be satisfied, in other words $U(D)=\emptyset$.

Finally, if $v_2(D)\ge5$, i.e., $D\equiv0\pmod{32}$, the above argument
shows that we must choose $d\equiv4\pmod8$ hence $d'\equiv0\pmod8$, and thus
the gcd condition is equivalent to $\gcd(d,d')=4$.

We have thus shown that there is a bijection between the cases $a=c$
and the elements $d\in U^+(D)$, and we have
$\tau=(-b+\sqrt{D})/(2a)=(-(|D|/d-d)+2\sqrt{D})/(|D|/d+d)$.
Now we check that $\tau_1=-1/(\tau+1)=-1/2+\sqrt{D}/(2d)$.
Thus $q_1=e^{2\pi i\tau_1}$ is a negative real, so from the product expansion
we deduce that $\eta^4(\tau_1)=re^{-2\pi i/12}$
for some positive real $r$. Thus by modularity
$\eta^4(-1/\tau_1)=\eta^4(\tau+1)=-\tau_1^2re^{-2\pi i/12}$, so
$\eta^4(\tau)=e^{-2\pi i/6}\eta^4(\tau+1)=e^{3i\pi/6}\tau_1^2r$.
Setting $r_1=r|\tau_1^2|$, we deduce that 
$$\eta^4(\tau)=e^{3i\pi/6}r_1\dfrac{\tau_1}{\ov{\tau_1}}=e^{3i\pi/6}r_1\dfrac{d-\sqrt{D}}{d+\sqrt{D}}\;,$$
for some other positive real number $r_1$.

\smallskip

The other case $a=b$ occurs for $D=a^2-4ac=-a(4c-a)$. Set $d=a$ and
$d'=4c-a=|D|/d$. The condition $a\le c$ is equivalent to $d\le\sqrt{|D|/3}$.
Since $c=(d+d')/4$, the integrality of $c$ together with
$\gcd(a,b,c)=\gcd(a,c)=1$ is equivalent to $\gcd(d+d',4d)=4$, and
since $4d-2(d+d')=2(d-d')$ this is equivalent to $\gcd(d+d',2(d-d'))=4$
hence to the same arithmetic conditions, so $d\in U^-(D)$. Here we have
$\tau=-1/2+\sqrt{D}/(2d)$, and as above, we deduce
from the product expansion that $\eta^4(\tau)=re^{-i\pi/6}$ for some positive
real $r$, proving the theorem since the product of the positive reals involved
is equal to $|P(D)|=1$ by the Chowla--Selberg formula.

Note that both cases can occur simultaneously only for $D=-3$, and one
checks that the result is still valid in that case.\fp

Recall that the choice of the standard set of representatives $F(D)$ is not
completely canonical: it corresponds to choosing the interior of the
standard fundamental domain of $\SL_2(\Z)$ together with the boundaries
the half-line $\Re(\tau)=-1/2$, $\Im(\tau)\ge\sqrt{3}/2$, and the small arc
$|\tau|=1$, $-1/2\le\Re(\tau)\le0$. But, if desired, we can modify $F(D)$ by
choosing as second boundary, instead of a small arc, the segment
$\Re(\tau)=1/2$ and $1/2\le\Im(\tau)\le\sqrt{3}/2$. This corresponds to
changing all occurrences of $a=c$,
i.e., the quadratic forms $(a,b,a)$, into $(2a-b,2a-b,a)$. If we call
$F'(D)$ this new set of representatives, and $P'(D)$ the corresponding
Chowla--Selberg quotient, by construction the set $U^+(D)$ disappears, and we
obtain the simpler formula $$P'(D)=e^{-(i\pi/6)|U(D)|}\;.$$

\begin{corollary}\label{cor1}
\begin{enumerate}
\item If $D\equiv0\pmod8$ we have $P(D)=1$.
\item Let $p$ be prime. If $D=-p$ with $p\equiv3\pmod4$ or $D=-4p$ with
$p\equiv1\pmod4$ we have $P(D)=e^{-2i\pi/12}$.
\end{enumerate}
\end{corollary}

\Proof (1) is clear since $U(D)=\emptyset$. For (2), if $D=-p$ with
$p\equiv3\pmod4$ we have $U(D)=\{1\}$ so $U^+(D)=\emptyset$ and $U^-(D)=\{1\}$
so the result follows. Similarly, if $D=-4p$ with $p\equiv1\pmod4$ we have
$U(D)=\{2\}$ and the result is also clear.\fp

\section{Application to Individual Values of $\eta(\tau)$}

Individual values of the modulus $|\eta(\tau)|$ have been given explicitly
starting from the pioneering work in \cite{vdP-Wil}, see also \cite{Cha-vdP}.
Here we simply give some immediate consequences of the above results for
$\eta(\tau)$ itself.

\begin{lemma}\label{lem1} If $\tau_1$ and $\tau_2$ are two CM points in the
same quadratic field $\Q(\sqrt{D_0})$ then $\eta(\tau_1)/\eta(\tau_2)$ is
an algebraic number.\end{lemma}

\Proof Since they are both in the same quadratic field, there exist
rational numbers $a\ne0$ and $b$ such that $\tau_2=a\tau_1+b$, so clearing
denominators, integers $A\ne0$, $B$, and $C\ne0$ such that
$C\tau_2=A\tau_1+B$. Now by a basic theorem of CM theory $\eta(C\tau_2)=\al\eta(\tau_2)$ for some algebraic number $\al$, and we also have
$$\eta(C\tau_2)=\eta(A\tau_1+B)=e^{2\pi iB/24}\eta(A\tau_1)=\be\eta(\tau_1)$$
for some algebraic number $\be$, so
$\eta(\tau_2)=(\be/\al)\eta(\tau_1)$.\fp

\begin{theorem}\label{thm3} Let $\tau$ be a CM point of discriminant $D$, and
as above write $D=D_0f^2$ with $D_0$ fundamental. There exists an algebraic
number $\al(\tau)$ such that
$$\Im(\tau)\eta(\tau)^4=\al(\tau)\dfrac{c(D_0,f)\CS(D_0)}{4\pi|D|^{1/2}}\;,$$
or, equivalently, if $\tau$ corresponds to $(a,b,c)$:
$$\eta(\tau)^4=a\al(\tau)\dfrac{c(D_0,f)\CS(D_0)}{2\pi|D|}\;.$$
\end{theorem}

\Proof By Theorems \ref{thm1} and \ref{thm2}, we have
$\prod_{\tau\in Q(D)}\eta(\tau)^4=\al_1(\CS(D_0)/\pi)^{h(D)}$
for some algebraic number $\al_1$. By the above lemma, all the factors
on the left-hand side are proportional up to an algebraic number, and since
$|Q(D)|=h(D)$ by definition, the result follows.\fp

The factors $\Im(\tau)$, $|D|^{1/2}$, and even $c(D_0,f)$ can of course be
removed without changing the validity of the theorem since they are
all algebraic, but in view of the Chowla--Selberg formula, the above
definition is the most natural, and in particular we know that the
product of the $\al(\tau)$ for $\tau\in F(D)$ has modulus $1$.

From this, we can trivially deduce the following general result which is
certainly classical:

\begin{proposition}\label{thm4} Let $f$ be a modular form or function of some
integral or half integral weight $k$ on some congruence subgroup of the
modular group, and assume that the coefficients of the Fourier expansion
of $f$ at infinity (or, for that matter, at any cusp) are algebraic.
For any CM point $\tau$, the number $f(\tau)/|f(\tau)|$ (hence its square
$f(\tau)/\ov{f(\tau)}$) is an algebraic number.\end{proposition}

\Proof Indeed, by CM theory we know that $f(\tau)/\eta(\tau)^{2k}$ is
an algebraic number, and by Theorem \ref{thm3} $\eta(\tau)/|\eta(\tau)|$
is algebraic, so the result follows.\fp

Note that since $(\eta(\tau)/|\eta(\tau)|)^2=\eta(\tau)/\ov{\eta(\tau)}=\eta(\tau)/\eta(-\ov{\tau})$ and $\tau$ and $-\ov{\tau}$ belong to the same
quadratic field, Lemma \ref{lem1} shows directly that
$\eta(\tau)/|\eta(\tau)|$ is algebraic without going through the explicit
computation given by Theorem \ref{thm2}. The same remark shows that
the proposition is trivially valid when the coefficients of the Fourier
expansion of $f$ are real since in that case $\ov{f(\tau)}=f(-\ov{\tau})$.

\medskip

To illustrate, we give a small table of values of the algebraic number
$\al(\tau)$ occurring in Theorem \ref{thm3}.

\bigskip

\centerline{
\begin{tabular}{|c||c|c|}
\hline
$D$ & $(a,b,c)$ & $\al(\tau)$ \\
\hline\hline
$-3$ & $(1,1,1)$ & $e^{-i\pi/6}$ \\
$-4$ & $(1,0,1)$ & $1$ \\
$-7$ & $(1,1,2)$ & $e^{-i\pi/6}$ \\
$-8$ & $(1,0,2)$ & $1$ \\
$-11$ & $(1,1,3)$ & $e^{-i\pi/6}$ \\
$-12$ & $(1,0,3)$ & $1$ \\
$-15$ & $(1,1,4)$ & $e^{-i\pi/6}((\sqrt{5}-1)/2)^{1/3}$ \\
      & $(2,1,2)$ & $e^{-i\pi/2}((\sqrt{-15}+1)/4)((\sqrt{5}+1)/2)^{1/3}$ \\
$-16$ & $(1,0,4)$ & $1$ \\
$-19$ & $(1,1,5)$ & $e^{-i\pi/6}$ \\
$-20$ & $(1,0,5)$ & $((\sqrt{5}-1)/2)^{1/2}$ \\
      & $(2,2,3)$ & $e^{-i\pi/6}((\sqrt{5}+1)/2)^{1/2}$\\
\hline\hline
\end{tabular}}

\bigskip

{\bf Remarks.}

\begin{enumerate}
\item By Corollary \ref{cor1}, for all odd fundamental discriminants $D$ of
  class number $1$, we have $\al(\tau)=e^{-i\pi/6}=(\sqrt{3}-\sqrt{-1})/2$,
  and for $8\mid D$ of class number $1$ we have $\al(\tau)=1$.
\item For values of $\tau$ with the above discriminants $D$ but not in the
  fundamental domain, we simply use the modularity of $\eta^4$
  to send $\tau$ to $F(D)$.
\item It is slightly surprising that the complex exponential which occurs
  is $e^{-i\pi/6}$ which is a $12$th root of unity, while since
  $\eta^4(\tau)=q^{1/6}(1+\cdots)=e^{2\pi i\tau/6}(1+\cdots)$,
  one would expect $6$th roots of unity instead, but in fact it should be
  interpreted as $e^{-i\pi/6}=i/\rho$, where $\rho=e^{2i\pi/3}$ is a cube
  root of unity.
\end{enumerate}

In view of this last remark, it seems useful to set the following definition:

\begin{definition} We define $\be(\tau)$ to be equal to $e^{mi\pi/6}\al(\tau)$
  for the smallest $m\ge0$ such that $\be(\tau)$ has minimal degree among
  the $12$ different values.\end{definition}

\section{Some Conjectures}

I have largely extended the above table (see below), and interestingly enough,
some precise conjectures have emerged, backed by a reasonable amount of data.

\begin{conjecture} Let $\tau=(-b+\sqrt{D})/(2a)$ with $D=b^2-4ac<0$, where
  as usual $a>0$ and $\gcd(a,b,c)=1$, and write $D=D_0f^2$ with $D_0$
  fundamental. As above, set
  $$\al(\tau)=\dfrac{\eta(\tau)^4}{a\cdot c(D_0,f)\CS(D_0)/(2\pi D)}$$
  and let $\be(\tau)$ as in the above definition.
  \begin{enumerate}
  \item The numbers $a\cdot\al(\tau)$ and $a/\al(\tau)$ are algebraic
    \emph{integers}.
  \item The norm of $\al(\tau)$ as an algebraic number (i.e., the product of
    its algebraic conjugates) is equal to $\pm1$ (since it is usually not an
    algebraic integer, this of course does not mean that it is a unit),
    and in fact perhaps always equal to $1$.
  \item Write $h(D)=2^vm$ with $m$ odd. The minimal polynomial of
    $a\cdot\be(\tau)$ is of the form $B(x^{jm})$ with $j=1$ if $v=0$,
    and $j\mid3\cdot2^{v-1}$ if $v\ge1$, and where the degree of $B$ is equal
    to $h(D)$ if $a=1$ and to $h(D)$ or $2h(D)$ if $a>1$.
  \item The number $\al(\tau)^{\lcm(6,h(D))}$ belongs to the ring class field
    of the quadratic order of discriminant $D$, of degree $2h(D)$, and
    when $b=0$ or $a=b$ it even belongs to its subfield of degree $h(D)$
    fixed by the automorphism $\sqrt{D}\mapsto-\sqrt{D}$.
  \end{enumerate}
\end{conjecture}

These conjectures are perhaps not difficult to prove using known results
on CM theory.

\medskip

{\bf Remarks.}

\begin{enumerate}
\item (1) and (2) mean that if $\al(\tau)$ is of degree $2d$,
up to a multiplicative constant its characteristic polynomial with integer
coefficients is of the form
$$a^dx^{2d}+c_{2d-1}a^{d-1}x^{2d-1}+\cdots+c_d x^d+c_{d-1}ax^{d-1}+\cdots+c_1a^{d-1}x+a^d\;.$$
If it has odd degree $2d-1$, it has a similar shape, but in all the
examples I have tested the polynomial is monic with constant term $-1$
(this may be true in general, but I do not have enough evidence to conjecture
it).
\item The exponent $\lcm(6,h(D))$ in (4) can be explained for at least
  two reasons. First, it may be reasonable to study
  $\al(\tau)^6$ instead of $\al(\tau)$, since this would correspond to
  $\Delta(\tau)$ instead of $\eta^4(\tau)$. Second, one could hope to
  get rid of this $6$ by considering $\be(\tau)$ instead of $\al(\tau)$,
  but for instance for $D=-15$, we really need an exponent divisible by
  $3$ in $\be(\tau)$ or $\al(\tau)$ to land in the required ring class field.
\item These conjectures imply that the degree of $\al(\tau)$ or of $\be(\tau)$
  is roughly proportional to $h(D)^2$, and in particular grows quite
  fast as $|D|\to\infty$. This is in contrast with the CM values of modular
  \emph{functions} such as $j(\tau)$, whose degrees are proportional to
  $h(D)$, not to its square.
\end{enumerate}

\medskip

These conjectures suggest the following algorithm to compute explicitly
$\al(\tau)$, hence to deduce an explicit expression for CM values of
$\eta(\tau)$:

\begin{enumerate}
\item Using the standard $\SL_2(\Z)$ transformations $\tau\mapsto\tau+1$
  and $\tau\mapsto-1/\tau$, we may assume that $\tau$ is in the standard
  fundamental domain (and to terminate, use
  $\eta(\tau+1)=e^{2\pi i/24}\eta(\tau)$ and
  $\eta(-1/\tau)=(\tau/i)^{1/2}\eta(\tau)$).
\item Write $\tau=(-b+\sqrt{D})/(2a)$ with $D=b^2-4ac<0$, $a>0$,
  and $\gcd(a,b,c)=1$, compute $h(D)$, write $D=D_0f^2$ with $D_0$
  fundamental, and compute $\al(\tau)$ and $\be(\tau)$ as defined above
  to sufficient accuracy. If the subsequent computations fail, recompute
  after increasing the accuracy, say by $50\%$ until they succeed.
\item Write $h(D)=2^vm$ with $m$ odd, and using an algebraic recognition
  program such as {\tt algdep} in {\tt Pari/GP}, check whether
  $A=(a\be(\tau))^{jm}$ is a root of a polynomial of degree $d$ for $j=1$ if
  $v=0$, or for $j$ a divisor of $3\cdot 2^{v-1}$ if $v\ge1$, with $d=h(D)$ if
  $a=1$ or $d=2h(D)$ if $a>1$.
  Factor the resulting polynomial over $\Z$, and check whether $A$ is an
  approximate root of a factor. If this is not the case, increase the
  accuracy and go back to the preceding step.
\item If $P(X)$ is the polynomial obtained in the previous step, the
  minimal polynomial of $a\be(\tau)$ will be a factor of $P(X^{jm})$,
  and from that it is immediate to obtain the minimal polynomial of
  $\al(\tau)$.
\end{enumerate}

\medskip

To illustrate, we have used this algorithm to construct a large table of
the degrees of the different expressions, where for notational simplicity
we set $d_{\al}=\deg(\al(\tau))$,
$d_{\be}=\deg(\be(\tau))$, $d_6=\deg(\al(\tau)^6)$, and
$d_{[6,h]}=\deg(\al(\tau)^{\lcm(6,h(D))})$. Conjecture (4) implies that the
last column is always a divisor of $2h$, and in addition a divisor of $h$
when $b=0$ or $b=a$ (it may also happen in other cases, the first
example being for $D=-76$). Note that we often have $d_6=d_{[6,h]}$,
but as soon as $8$ or some prime $p\ge5$ divides $h(D)$, we have
$d_6>d_{[6,h]}$, and this happens for almost all $D$.

The indicated degrees $d$ are either coded as $e\cdot f$, meaning
that it is a polynomial of degree $e$ in $x^f$ (hence of degree $ef$),
or simply as $d$ if $f=1$.

\vfill\eject

\bigskip

\centerline{
\begin{tabular}{|c||c||c|c|c|c|}
\hline
$D$ & $(a,b,c)$ & $d_{\al}$ & $d_{\be}$ & $d_6$ & $d_{[6,h]}$ \\
\hline\hline
$-3$ & $(1,1,1)$ & 2$\cdot$2 & 1 & 1 & 1 \\
$-4$ & $(1,0,1)$ & 1 & 1 & 1 & 1 \\
$-7$ & $(1,1,2)$ & 2$\cdot$2 & 1 & 1 & 1 \\
$-8$ & $(1,0,2)$ & 1 & 1 & 1 & 1 \\
$-11$ & $(1,1,3)$ & 2$\cdot$2 & 1 & 1 & 1 \\
$-12$ & $(1,0,3)$ & 1 & 1 & 1 & 1 \\
$-15$ & $(1,1,4)$ & 2$\cdot$6 & 2$\cdot$3 & 2 & 2 \\
      & $(2,1,2)$ & 4$\cdot$6 & 4$\cdot$3 & 4 & 4 \\
$-16$ & $(1,0,4)$ & 1 & 1 & 1 & 1 \\
$-19$ & $(1,1,5)$ & 2$\cdot$2 & 1 & 1 & 1 \\
$-20$ & $(1,0,5)$ & 2$\cdot$2 & 2$\cdot$2 & 2 & 2 \\
      & $(2,2,3)$ & 4$\cdot$2 & 2$\cdot$2 & 2 & 2 \\
$-23$ & $(1,1,6)$ & 3$\cdot$6 & 3$\cdot$3 & 3 & 3 \\
      & $(2,\pm1,3)$ & 6$\cdot$6 & 6$\cdot$3 & 6 & 6 \\
$-24$ & $(1,0,6)$ & 2$\cdot$3 & 2$\cdot$3 & 2 & 2 \\
      & $(2,0,3)$ & 2$\cdot$3 & 2$\cdot$3 & 2 & 2 \\
$-27$ & $(1,1,7)$ & 2$\cdot$2 & 1 & 1 & 1 \\
$-28$ & $(1,0,7)$ & 1 & 1 & 1 & 1 \\
$-31$ & $(1,1,8)$ & 3$\cdot$6 & 3$\cdot$3 & 3 & 3 \\
      & $(2,\pm1,4)$ & 6$\cdot$6 & 6$\cdot$3 & 6 & 6 \\
$-32$ & $(1,0,8)$ & 2$\cdot$2 & 2$\cdot$2 & 2 & 2 \\
      & $(3,2,3)$ & 4$\cdot$2 & 4$\cdot$2 & 4 & 4 \\
$-35$ & $(1,1,9)$ & 4$\cdot$2 & 2 & 2 & 2 \\
      & $(3,1,3)$ & 4$\cdot$2 & 4 & 4 & 4 \\
$-36$ & $(1,0,9)$ & 2$\cdot$3 & 2$\cdot$3 & 2 & 2 \\
      & $(2,2,5)$ & 2$\cdot$6 & 2$\cdot$3 & 2 & 2 \\
$-39$ & $(1,1,10)$ & 4$\cdot$6 & 4$\cdot$6 & 4 & 4 \\
      & $(2,\pm1,5)$ & 8$\cdot$6 & 8$\cdot$6 & 8 & 8 \\
      & $(3,3,4)$ & 4$\cdot$6 & 4$\cdot$6 & 4 & 4 \\
$-40$ & $(1,0,10)$ & 2 & 2 & 2 & 2 \\
      & $(2,0,5)$ & 2 & 2 & 2 & 2 \\
$-43$ & $(1,1,11)$ & 2$\cdot$2 & 1 & 1 & 1 \\
$-44$ & $(1,0,11)$ & 3$\cdot$3 & 3$\cdot$3 & 3 & 3 \\
      & $(3,\pm2,4)$ & 6$\cdot$3 & 6$\cdot$3 & 6 & 6 \\
$-47$ & $(1,1,12)$ & 10$\cdot$10 & 5$\cdot$5 & 5$\cdot$5 & 5 \\
      & $(2,\pm1,6)$ & 20$\cdot$10 & 10$\cdot$5 & 10$\cdot$5 & 10 \\
      & $(3,\pm1,4)$ & 20$\cdot$10 & 10$\cdot$5 & 10$\cdot$5 & 10 \\
$-48$ & $(1,0,12)$ & 2$\cdot$2 & 2$\cdot$2 & 2 & 2 \\
      & $(3,0,4)$ & 2$\cdot$2 & 2$\cdot$2 & 2 & 2 \\
\hline\hline
\end{tabular}}
\vfill\eject

\bigskip

\centerline{
\begin{tabular}{|c||c||c|c|c|c|}
\hline
$D$ & $(a,b,c)$ & $d_{\al}$ & $d_{\be}$ & $d_6$ & $d_{[6,h]}$ \\
\hline\hline
$-51$ & $(1,1,13)$ & 2$\cdot$6 & 2$\cdot$3 & 2 & 2 \\
      & $(3,3,5)$ & 2$\cdot$6 & 2$\cdot$3 & 2 & 2 \\
$-52$ & $(1,0,13)$ & 2$\cdot$2 & 2$\cdot$2 & 2 & 2 \\
      & $(2,2,7)$ & 4$\cdot$2 & 2$\cdot$2 & 2 & 2 \\
$-55$ & $(1,1,14)$ & 8$\cdot$2 & 4$\cdot$2 & 4 & 4 \\
      & $(2,\pm1,7)$ & 16$\cdot$2 & 8$\cdot$2 & 8 & 8 \\
      & $(4,3,4)$ & 8$\cdot$2 & 8$\cdot$2 & 8 & 8 \\
$-56$ & $(1,0,14)$ & 4$\cdot$2 & 4$\cdot$2 & 4 & 4 \\
      & $(2,0,7)$ & 4$\cdot$2 & 4$\cdot$2 & 4 & 4 \\
      & $(3,\pm2,5)$ & 16$\cdot$2 & 8$\cdot$2 & 8 & 8 \\
$-59$ & $(1,1,15)$ & 3$\cdot$6 & 3$\cdot$3 & 3 & 3 \\
      & $(3,\pm1,5)$ & 6$\cdot$6 & 6$\cdot$3 & 6 & 6 \\
$-60$ & $(1,0,15)$ & 2$\cdot$3 & 2$\cdot$3 & 2 & 2 \\
      & $(3,0,5)$ & 2$\cdot$3 & 2$\cdot$3 & 2 & 2 \\
$-63$ & $(1,1,16)$ & 4$\cdot$6 & 4$\cdot$6 & 4 & 4 \\
      & $(2,\pm1,8)$ & 8$\cdot$3 & 8$\cdot$3 & 8 & 8 \\
      & $(4,1,4)$ & 8$\cdot$3 & 8$\cdot$3 & 8 & 8 \\
$-64$ & $(1,0,16)$ & 2 & 2 & 2 & 2 \\
      & $(4,4,5)$ & 4$\cdot$2 & 2 & 2 & 2 \\
$-67$ & $(1,1,17)$ & 2$\cdot$2 & 1 & 1 & 1 \\
$-68$ & $(1,0,17)$ & 4$\cdot$2 & 4$\cdot$2 & 4 & 4 \\
      & $(2,2,9)$ & 8$\cdot$2 & 4$\cdot$2 & 4 & 4 \\
      & $(3,\pm2,6)$ & 8$\cdot$2 & 8$\cdot$2 & 8 & 8 \\
$-71$ & $(1,1,18)$ & 14$\cdot$14 & 7$\cdot$7 & 7$\cdot$7 & 7 \\
      & $(2,\pm1,9)$ & 28$\cdot$14 & 14$\cdot$7 & 14$\cdot$7 & 14 \\
      & $(3,\pm1,6)$ & 14$\cdot$14 & 14$\cdot$7 & 14$\cdot$7 & 14 \\
      & $(4,\pm3,5)$ & 14$\cdot$14 & 14$\cdot$7 & 14$\cdot$7 & 14 \\
$-72$ & $(1,0,18)$ & 2$\cdot$3 & 2$\cdot$3 & 2 & 2 \\
      & $(2,0,9)$ & 2$\cdot$3 & 2$\cdot$3 & 2 & 2 \\
$-75$ & $(1,1,19)$ & 4$\cdot$2 & 2 & 2 & 2 \\
      & $(3,3,7)$ & 4$\cdot$2 & 2 & 2 & 2 \\
$-76$ & $(1,0,19)$ & 3$\cdot$3 & 3$\cdot$3 & 3 & 3 \\
      & $(4,\pm2,5)$ & 3$\cdot$3 & 3$\cdot$3 & 3 & 3 \\
$-79$ & $(1,1,20)$ & 10$\cdot$10 & 5$\cdot$5 & 5$\cdot$5 & 5 \\
      & $(2,\pm1,10)$ & 20$\cdot$10 & 10$\cdot$5 & 10$\cdot$5 & 10 \\
      & $(4,\pm1,5)$ & 20$\cdot$10 & 10$\cdot$5 & 10$\cdot$5 & 10 \\
$-80$ & $(1,0,20)$ & 4 & 4 & 4 & 4 \\
      & $(3,\pm2,7)$ & 16 & 8 & 8 & 8 \\
      & $(4,0,5)$ & 4 & 4 & 4 & 4 \\
\hline\hline
\end{tabular}}
\vfill\eject

\bigskip

\centerline{
\begin{tabular}{|c||c||c|c|c|c|}
\hline
$D$ & $(a,b,c)$ & $d_{\al}$ & $d_{\be}$ & $d_6$ & $d_{[6,h]}$ \\
\hline\hline
$-83$ & $(1,1,21)$ & 3$\cdot$6 & 3$\cdot$3 & 3 & 3 \\
      & $(3,\pm1,7)$ & 6$\cdot$6 & 6$\cdot$3 & 6 & 6 \\
$-84$ & $(1,0,21)$ & 4$\cdot$6 & 4$\cdot$6 & 4 & 4 \\
      & $(2,2,11)$ & 4$\cdot$6 & 4$\cdot$6 & 4 & 4 \\
      & $(3,0,7)$ & 4$\cdot$6 & 4$\cdot$6 & 4 & 4 \\
      & $(5,4,5)$ & 8$\cdot$6 & 8$\cdot$6 & 8 & 8 \\
$-87$ & $(1,1,22)$ & 6$\cdot$6 & 6$\cdot$3 & 6 & 6 \\
      & $(2,\pm1,11)$ & 12$\cdot$6 & 12$\cdot$3 & 12 & 12 \\
      & $(3,3,8)$ & 6$\cdot$6 & 6$\cdot$3 & 6 & 6 \\
      & $(4,\pm3,6)$ & 12$\cdot$6 & 12$\cdot$3 & 12 & 12 \\
$-88$ & $(1,0,22)$ & 2 & 2 & 2 & 2 \\
      & $(2,0,11)$ & 2 & 2 & 2 & 2 \\
$-91$ & $(1,1,23)$ & 4$\cdot$2 & 2 & 2 & 2 \\
      & $(5,3,5)$ & 4$\cdot$2 & 4 & 4 & 4 \\
$-92$ & $(1,0,23)$ & 3$\cdot$3 & 3$\cdot$3 & 3 & 3 \\
      & $(3,\pm2,8)$ & 6$\cdot$3 & 6$\cdot$3 & 6 & 6 \\
$-95$ & $(1,1,24)$ & 16$\cdot$4 & 8$\cdot$4 & 8$\cdot$2 & 8 \\
      & $(2,\pm1,12)$ & 32$\cdot$4 & 16$\cdot$4 & 16$\cdot$2 & 16 \\
      & $(3,\pm1,8)$ & 32$\cdot$4 & 16$\cdot$4 & 16$\cdot$2 & 16 \\
      & $(4,\pm1,6)$ & 32$\cdot$4 & 16$\cdot$4 & 16$\cdot$2 & 16 \\
      & $(5,5,6)$ & 16$\cdot$4 & 8$\cdot$4 & 8$\cdot$2 & 8 \\
$-96$ & $(1,0,24)$ & 4$\cdot$6 & 4$\cdot$6 & 4 & 4 \\
      & $(3,0,8)$ & 4$\cdot$6 & 4$\cdot$6 & 4 & 4 \\
      & $(4,4,7)$ & 4$\cdot$6 & 4$\cdot$6 & 4 & 4 \\
      & $(5,2,5)$ & 8$\cdot$6 & 8$\cdot$6 & 8 & 8 \\
$-99$ & $(1,1,25)$ & 2$\cdot$6 & 2$\cdot$3 & 2 & 2 \\
      & $(5,1,5)$ & 4$\cdot$6 & 4$\cdot$3 & 4 & 4 \\
$-100$ & $(1,0,25)$ & 2 & 2 & 2 & 2 \\
      & $(2,2,13)$ & 4$\cdot$2 & 2 & 2 & 2 \\
$-103$ & $(1,1,26)$ & 10$\cdot$10 & 5$\cdot$5 & 5$\cdot$5 & 5 \\
      & $(2,\pm1,13)$ & 20$\cdot$10 & 10$\cdot$5 & 10$\cdot$5 & 10 \\
      & $(4,\pm3,7)$ & 10$\cdot$10 & 10$\cdot$5 & 10$\cdot$5 & 10 \\
$-104$ & $(1,0,26)$ & 6$\cdot$3 & 6$\cdot$3 & 6 & 6 \\
      & $(2,0,13)$ & 6$\cdot$3 & 6$\cdot$3 & 6 & 6 \\
      & $(3,\pm2,9)$ & 12$\cdot$3 & 12$\cdot$3 & 12 & 12 \\
      & $(5,\pm4,6)$ & 12$\cdot$3 & 12$\cdot$3 & 12 & 12 \\
$-107$ & $(1,1,27)$ & 3$\cdot$6 & 3$\cdot$3 & 3 & 3 \\
      & $(3,\pm1,9)$ & 6$\cdot$6 & 6$\cdot$3 & 6 & 6 \\
$-108$ & $(1,0,27)$ & 3$\cdot$3 & 3$\cdot$3 & 3 & 3 \\
      & $(4,\pm2,7)$ & 3$\cdot$3 & 3$\cdot$3 & 3 & 3 \\
\hline\hline
\end{tabular}}
\vfill\eject

\bigskip

\centerline{
\begin{tabular}{|c||c||c|c|c|c|}
\hline
$D$ & $(a,b,c)$ & $d_{\al}$ & $d_{\be}$ & $d_6$ & $d_{[6,h]}$ \\
\hline\hline
$-111$ & $(1,1,28)$ & 8$\cdot$12 & 8$\cdot$12 & 8$\cdot$2 & 8 \\
      & $(2,\pm1,14)$ & 16$\cdot$12 & 16$\cdot$12 & 16$\cdot$2 & 16 \\
      & $(3,3,10)$ & 8$\cdot$12 & 8$\cdot$12 & 8$\cdot$2 & 8 \\
      & $(4,\pm1,7)$ & 16$\cdot$12 & 16$\cdot$12 & 16$\cdot$2 & 16 \\
      & $(5,\pm3,6)$ & 16$\cdot$12 & 16$\cdot$12 & 16$\cdot$2 & 16 \\
$-112$ & $(1,0,28)$ & 2$\cdot$2 & 2$\cdot$2 & 2 & 2 \\
      & $(4,0,7)$ & 2$\cdot$2 & 2$\cdot$2 & 2 & 2 \\
$-115$ & $(1,1,29)$ & 4$\cdot$2 & 2 & 2 & 2 \\
      & $(5,5,7)$ & 4$\cdot$2 & 2 & 2 & 2 \\
$-116$ & $(1,0,29)$ & 6$\cdot$6 & 6$\cdot$6 & 6 & 6 \\
      & $(2,2,15)$ & 6$\cdot$6 & 6$\cdot$6 & 6 & 6 \\
      & $(3,\pm2,10)$ & 12$\cdot$6 & 12$\cdot$6 & 12 & 12 \\
      & $(5,\pm2,6)$ & 12$\cdot$6 & 12$\cdot$6 & 12 & 12 \\
$-119$ & $(1,1,30)$ & 20$\cdot$10 & 10$\cdot$5 & 10$\cdot$5 & 10 \\
      & $(2,\pm1,15)$ & 40$\cdot$10 & 20$\cdot$5 & 20$\cdot$5 & 20 \\
      & $(3,\pm1,10)$ & 40$\cdot$10 & 20$\cdot$5 & 20$\cdot$5 & 20 \\
      & $(4,\pm3,8)$ & 20$\cdot$10 & 20$\cdot$5 & 20$\cdot$5 & 20 \\
      & $(5,\pm1,6)$ & 40$\cdot$10 & 20$\cdot$5 & 20$\cdot$5 & 20 \\
      & $(6,5,6)$ & 20$\cdot$10 & 20$\cdot$5 & 20$\cdot$5 & 20 \\
$-120$ & $(1,0,30)$ & 4$\cdot$3 & 4$\cdot$3 & 4 & 4 \\
      & $(2,0,15)$ & 4$\cdot$3 & 4$\cdot$3 & 4 & 4 \\
      & $(3,0,10)$ & 4$\cdot$3 & 4$\cdot$3 & 4 & 4 \\
      & $(5,0,6)$ & 4$\cdot$3 & 4$\cdot$3 & 4 & 4 \\
$-123$ & $(1,1,31)$ & 2$\cdot$6 & 2$\cdot$3 & 2 & 2 \\
      & $(3,3,11)$ & 2$\cdot$6 & 2$\cdot$3 & 2 & 2 \\
$-124$ & $(1,0,31)$ & 3$\cdot$3 & 3$\cdot$3 & 3 & 3 \\
      & $(5,\pm4,7)$ & 6$\cdot$3 & 6$\cdot$3 & 6 & 6 \\
$-127$ & $(1,1,32)$ & 10$\cdot$10 & 5$\cdot$5 & 5$\cdot$5 & 5 \\
      & $(2,\pm1,16)$ & 20$\cdot$10 & 10$\cdot$5 & 10$\cdot$5 & 10 \\
      & $(4,\pm1,8)$ & 20$\cdot$10 & 10$\cdot$5 & 10$\cdot$5 & 10 \\
$-128$ & $(1,0,32)$ & 4$\cdot$2 & 4$\cdot$2 & 4 & 4 \\
      & $(3,\pm2,11)$ & 16$\cdot$2 & 8$\cdot$2 & 8 & 8 \\
      & $(4,4,9)$ & 8$\cdot$2 & 4$\cdot$2 & 4 & 4 \\
$-131$ & $(1,1,33)$ & 10$\cdot$10 & 5$\cdot$5 & 5$\cdot$5 & 5 \\
      & $(3,\pm1,11)$ & 20$\cdot$10 & 10$\cdot$5 & 10$\cdot$5 & 10 \\
      & $(5,\pm3,7)$ & 10$\cdot$10 & 10$\cdot$5 & 10$\cdot$5 & 10 \\
$-132$ & $(1,0,33)$ & 4$\cdot$6 & 4$\cdot$6 & 4 & 4 \\
      & $(2,2,17)$ & 4$\cdot$6 & 4$\cdot$6 & 4 & 4 \\
      & $(3,0,11)$ & 4$\cdot$6 & 4$\cdot$6 & 4 & 4 \\
      & $(6,6,7)$ & 4$\cdot$6 & 4$\cdot$6 & 4 & 4 \\
\hline\hline
\end{tabular}}
\vfill\eject

\bigskip

\centerline{
\begin{tabular}{|c||c||c|c|c|c|}
\hline
$D$ & $(a,b,c)$ & $d_{\al}$ & $d_{\be}$ & $d_6$ & $d_{[6,h]}$ \\
\hline\hline
$-135$ & $(1,1,34)$ & 6$\cdot$6 & 6$\cdot$3 & 6 & 6 \\
      & $(2,\pm1,17)$ & 12$\cdot$6 & 12$\cdot$3 & 12 & 12 \\
      & $(4,\pm3,9)$ & 12$\cdot$6 & 12$\cdot$3 & 12 & 12 \\
      & $(5,5,8)$ & 6$\cdot$6 & 6$\cdot$3 & 6 & 6 \\
$-136$ & $(1,0,34)$ & 4$\cdot$2 & 4$\cdot$2 & 4 & 4 \\
      & $(2,0,17)$ & 4$\cdot$2 & 4$\cdot$2 & 4 & 4 \\
      & $(5,\pm2,7)$ & 16$\cdot$2 & 8$\cdot$2 & 8 & 8 \\
$-139$ & $(1,1,35)$ & 3$\cdot$6 & 3$\cdot$3 & 3 & 3 \\
      & $(5,\pm1,7)$ & 6$\cdot$6 & 6$\cdot$3 & 6 & 6 \\
$-140$ & $(1,0,35)$ & 6$\cdot$3 & 6$\cdot$3 & 6 & 6 \\
      & $(3,\pm2,12)$ & 12$\cdot$3 & 12$\cdot$3 & 12 & 12 \\
      & $(4,\pm2,9)$ & 6$\cdot$3 & 6$\cdot$3 & 6 & 6 \\
      & $(5,0,7)$ & 6$\cdot$3 & 6$\cdot$3 & 6 & 6 \\
$-143$ & $(1,1,36)$ & 20$\cdot$10 & 10$\cdot$5 & 10$\cdot$5 & 10 \\
      & $(2,\pm1,18)$ & 40$\cdot$10 & 20$\cdot$5 & 20$\cdot$5 & 20 \\
      & $(3,\pm1,12)$ & 20$\cdot$10 & 20$\cdot$5 & 20$\cdot$5 & 20 \\
      & $(4,\pm1,9)$ & 40$\cdot$10 & 20$\cdot$5 & 20$\cdot$5 & 20 \\
      & $(6,1,6)$ & 20$\cdot$10 & 20$\cdot$5 & 20$\cdot$5 & 20 \\
      & $(6,\pm5,7)$ & 40$\cdot$10 & 20$\cdot$5 & 20$\cdot$5 & 20 \\
$-144$ & $(1,0,36)$ & 4$\cdot$3 & 4$\cdot$3 & 4 & 4 \\
      & $(4,0,9)$ & 4$\cdot$3 & 4$\cdot$3 & 4 & 4 \\
      & $(5,\pm4,8)$ & 8$\cdot$3 & 8$\cdot$3 & 8 & 8 \\
$-147$ & $(1,1,37)$ & 4$\cdot$2 & 2 & 2 & 2 \\
      & $(3,3,13)$ & 4$\cdot$2 & 2 & 2 & 2 \\
$-148$ & $(1,0,37)$ & 2$\cdot$2 & 2$\cdot$2 & 2 & 2 \\
      & $(2,2,19)$ & 4$\cdot$2 & 2$\cdot$2 & 2 & 2 \\
$-151$ & $(1,1,38)$ & 14$\cdot$14 & 7$\cdot$7 & 7$\cdot$7 & 7 \\
      & $(2,\pm1,19)$ & 28$\cdot$14 & 14$\cdot$7 & 14$\cdot$7 & 14 \\
      & $(4,\pm3,10)$ & 14$\cdot$14 & 14$\cdot$7 & 14$\cdot$7 & 14 \\
      & $(5,\pm3,8)$ & 14$\cdot$14 & 14$\cdot$7 & 14$\cdot$7 & 14 \\
$-152$ & $(1,0,38)$ & 6$\cdot$3 & 6$\cdot$3 & 6 & 6 \\
      & $(2,0,19)$ & 6$\cdot$3 & 6$\cdot$3 & 6 & 6 \\
      & $(3,\pm2,13)$ & 12$\cdot$3 & 12$\cdot$3 & 12 & 12 \\
      & $(6,\pm4,7)$ & 12$\cdot$3 & 12$\cdot$3 & 12 & 12 \\
$-155$ & $(1,1,39)$ & 8$\cdot$2 & 4$\cdot$2 & 4 & 4 \\
      & $(3,\pm1,13)$ & 16$\cdot$2 & 8$\cdot$2 & 8 & 8 \\
      & $(5,5,9)$ & 8$\cdot$2 & 4$\cdot$2 & 4 & 4 \\
$-156$ & $(1,0,39)$ & 4$\cdot$6 & 4$\cdot$6 & 4 & 4 \\
      & $(3,0,13)$ & 4$\cdot$6 & 4$\cdot$6 & 4 & 4 \\
      & $(5,\pm2,8)$ & 8$\cdot$6 & 8$\cdot$6 & 8 & 8 \\
\hline\hline
\end{tabular}}
\vfill\eject

\bigskip

\centerline{
\begin{tabular}{|c||c||c|c|c|c|}
\hline
$D$ & $(a,b,c)$ & $d_{\al}$ & $d_{\be}$ & $d_6$ & $d_{[6,h]}$ \\
\hline\hline
$-159$ & $(1,1,40)$ & 10$\cdot$30 & 10$\cdot$15 & 10$\cdot$5 & 10 \\
      & $(2,\pm1,20)$ & 20$\cdot$30 & 20$\cdot$15 & 20$\cdot$5 & 20 \\
      & $(3,3,14)$ & 10$\cdot$30 & 10$\cdot$15 & 10$\cdot$5 & 10 \\
      & $(4,\pm1,10)$ & 20$\cdot$30 & 20$\cdot$15 & 20$\cdot$5 & 20 \\
      & $(5,\pm1,8)$ & 20$\cdot$30 & 20$\cdot$15 & 20$\cdot$5 & 20 \\
      & $(6,\pm3,7)$ & 20$\cdot$30 & 20$\cdot$15 & 20$\cdot$5 & 20 \\
$-160$ & $(1,0,40)$ & 4$\cdot$2 & 4$\cdot$2 & 4 & 4 \\
      & $(4,4,11)$ & 8$\cdot$2 & 4$\cdot$2 & 4 & 4 \\
      & $(5,0,8)$ & 4$\cdot$2 & 4$\cdot$2 & 4 & 4 \\
      & $(7,6,7)$ & 8$\cdot$2 & 8$\cdot$2 & 8 & 8 \\
$-163$ & $(1,1,41)$ & 2$\cdot$2 & 1 & 1 & 1 \\
$-164$ & $(1,0,41)$ & 8$\cdot$4 & 8$\cdot$4 & 8$\cdot$2 & 8 \\
      & $(2,2,21)$ & 16$\cdot$4 & 8$\cdot$4 & 8$\cdot$2 & 8 \\
      & $(3,\pm2,14)$ & 32$\cdot$4 & 16$\cdot$4 & 16$\cdot$2 & 16 \\
      & $(5,\pm4,9)$ & 32$\cdot$4 & 16$\cdot$4 & 16$\cdot$2 & 16 \\
      & $(6,\pm2,7)$ & 32$\cdot$4 & 16$\cdot$4 & 16$\cdot$2 & 16 \\
$-167$ & $(1,1,42)$ & 22$\cdot$22 & 11$\cdot$11 & 11$\cdot$11 & 11 \\
      & $(2,\pm1,21)$ & 44$\cdot$22 & 22$\cdot$11 & 22$\cdot$11 & 22 \\
      & $(3,\pm1,14)$ & 44$\cdot$22 & 22$\cdot$11 & 22$\cdot$11 & 22 \\
      & $(4,\pm3,11)$ & 22$\cdot$22 & 22$\cdot$11 & 22$\cdot$11 & 22 \\
      & $(6,\pm1,7)$ & 44$\cdot$22 & 22$\cdot$11 & 22$\cdot$11 & 22 \\
      & $(6,\pm5,8)$ & 44$\cdot$22 & 22$\cdot$11 & 22$\cdot$11 & 22 \\
$-168$ & $(1,0,42)$ & 4$\cdot$6 & 4$\cdot$6 & 4 & 4 \\
      & $(2,0,21)$ & 4$\cdot$6 & 4$\cdot$6 & 4 & 4 \\
      & $(3,0,14)$ & 4$\cdot$6 & 4$\cdot$6 & 4 & 4 \\
      & $(6,0,7)$ & 4$\cdot$6 & 4$\cdot$6 & 4 & 4 \\
$-171$ & $(1,1,43)$ & 4$\cdot$6 & 4$\cdot$6 & 4 & 4 \\
      & $(5,\pm3,9)$ & 8$\cdot$3 & 8$\cdot$3 & 8 & 8 \\
      & $(7,5,7)$ & 8$\cdot$3 & 8$\cdot$3 & 8 & 8 \\
$-172$ & $(1,0,43)$ & 3$\cdot$3 & 3$\cdot$3 & 3 & 3 \\
      & $(4,\pm2,11)$ & 3$\cdot$3 & 3$\cdot$3 & 3 & 3 \\
$-175$ & $(1,1,44)$ & 6$\cdot$6 & 6$\cdot$3 & 6 & 6 \\
      & $(2,\pm1,22)$ & 12$\cdot$6 & 12$\cdot$3 & 12 & 12 \\
      & $(4,\pm1,11)$ & 12$\cdot$6 & 12$\cdot$3 & 12 & 12 \\
      & $(7,7,8)$ & 6$\cdot$6 & 6$\cdot$3 & 6 & 6 \\
$-176$ & $(1,0,44)$ & 6$\cdot$2 & 6$\cdot$2 & 6 & 6 \\
      & $(3,\pm2,15)$ & 12$\cdot$2 & 12$\cdot$2 & 12 & 12 \\
      & $(4,0,11)$ & 6$\cdot$2 & 6$\cdot$2 & 6 & 6 \\
      & $(5,\pm2,9)$ & 24$\cdot$2 & 12$\cdot$2 & 12 & 12 \\
$-179$ & $(1,1,45)$ & 10$\cdot$10 & 5$\cdot$5 & 5$\cdot$5 & 5 \\
      & $(3,\pm1,15)$ & 10$\cdot$10 & 10$\cdot$5 & 10$\cdot$5 & 10 \\
      & $(5,\pm1,9)$ & 20$\cdot$10 & 10$\cdot$5 & 10$\cdot$5 & 10 \\
$-180$ & $(1,0,45)$ & 4$\cdot$6 & 4$\cdot$6 & 4 & 4 \\
      & $(2,2,23)$ & 4$\cdot$6 & 4$\cdot$6 & 4 & 4 \\
      & $(5,0,9)$ & 4$\cdot$6 & 4$\cdot$6 & 4 & 4 \\
      & $(7,4,7)$ & 8$\cdot$6 & 8$\cdot$6 & 8 & 8 \\
\hline\hline
\end{tabular}}
\vfill\eject

\end{document}